\documentclass[a4paper]{amsart}
\usepackage{amsmath}
\usepackage{amsfonts}
\usepackage{amsthm}
\usepackage{tikz-cd}
\usepackage{amssymb}
\usepackage{hyperref}

\newcommand{\Q}{\mathbb{Q}}
\newcommand{\R}{\mathbb{R}}
\newcommand{\N}{\mathbb{N}}
\newcommand{\Z}{\mathbb{Z}}
\renewcommand{\P}{\mathbb{P}}
\newcommand{\F}{\mathbb{F}}
\renewcommand{\O}{\mathcal{O}}
\newcommand{\ad}{\mathbf{A}}
\newcommand{\pl}{\Omega}
\newcommand{\A}{\mathbb{A}}
\newcommand{\kbar}{{\bar{k}}}
\newcommand{\divs}{\mid}
\newcommand{\ndivs}{\nmid}
\newcommand{\Mat}{\mathbf{M}}
\newcommand{\HH}{\mathbb{H}}
\newcommand{\op}{\mathrm{op}}
\renewcommand{\H}{\mathrm{H}}
\newcommand{\al}{\mathcal{A}}
\newcommand{\m}{\mathfrak{m}}
\newcommand{\et}{\textrm{\'et}}
\newcommand{\Gm}{\mathbb{G}_{\mathrm{m}}}
\newcommand{\Xbar}{\bar{X}}
\newcommand{\nr}{\mathrm{nr}}
\newcommand{\Y}{\mathcal{Y}}
\newcommand{\np}{(p')}

\DeclareMathOperator{\Br}{Br}
\DeclareMathOperator{\Pic}{Pic}
\DeclareMathOperator{\End}{End}

\DeclareMathOperator{\im}{im}
\DeclareMathOperator{\Spec}{Spec}
\DeclareMathOperator{\Gal}{Gal}

\DeclareMathOperator{\Hom}{Hom}

\DeclareMathOperator{\inv}{inv}
\DeclareMathOperator{\Norm}{N}

\newtheorem{theorem}{Theorem}[section]

\theoremstyle{remark}
\newtheorem{remark}[theorem]{Remark}
\theoremstyle{definition}
\newtheorem{definition}[theorem]{Definition}
\newtheorem{example}[theorem]{Example}
\numberwithin{equation}{section}

\title{Rational points in families of varieties}
\author{Martin Bright}

\begin{document}

\maketitle

\tableofcontents

These are notes from lectures given at the Autumn School on Algebraic and Arithmetic Geometry at the Johannes Gutenberg-Universit\"at Mainz in October 2017.
I thank for organisers of the school, Dino Festi, Ariyan Javanpeykar and Davide Veniani, for their hospitality.

\section{The Hasse principle}

Let $X$ be a projective variety over the rational numbers $\Q$.  In order to study the rational points of $X$, we can embed $\Q$ in each of its completions, which are the real numbers $\R$ and the $p$-adic numbers $\Q_p$ for every prime $p$.  This gives an embedding of $X(\Q)$ into $X(\R)$ and into each of the $X(\Q_p)$.  It follows that a necessary condition for $X(\Q)$ to be non-empty is that $X(\R)$ and each $X(\Q_p)$ be non-empty.

This necessary condition is useful because it is a condition that can be checked algorithmically.

\subsection{Testing local solubility}

First consider testing, for a single prime $p$, whether $X$ has any $\Q_p$-points.  
For any fixed $r$, the finite set $X(\Z/p^r\Z)$ can be enumerated.
Hensel's Lemma states that, as soon as we can find a point of $X(\Z/p^r\Z)$ at which the derivatives of the defining polynomials of $X$ satisfy certain conditions, then there exists a point of $X(\Q_p)$.  
It is not immediately clear that this leads to a finite procedure for checking solubility, since we might end up looking at $X(\Z/p^r\Z)$ for arbitrarily large $r$.  However, it is easy to prove that failure to find a point satisfying the conditions of Hensel's Lemma leads to a sequence converging to a singular point of $X(\Q_p)$.  If $X$ is smooth, this gives a contradiction.  If $X$ is not smooth, then we can inductively look first for points on the singular subvariety of $X$, and then for smooth points on $X$.
In fact, it is possible to give an \emph{a priori} bound for the largest $r$ which is needed: see~\cite{Greenberg:IHES-1966}.

\begin{remark}
Sophisticated geometers will have noticed that, to talk about points of $X(\Z/p\Z)$, we need to first fix a model of $X$ over $\Z$, or at least over $\Z_p$.
However, in this section we are thinking of $X$ as being given explicitly in projective space by equations, so it already comes with a choice of model (given by taking the Zariski closure of $X \subset \P^n_\Q$ in $\P^n_\Z$).  We can therefore talk about points of $X$ over any ring without worrying.
\end{remark}

We will make a lot of use of the following well-known special case of Hensel's Lemma: if $X$ has a smooth point over $\F_p = \Z/p\Z$, then $X$ has a $\Z_p$-point.  More generally, we have the following.

\begin{theorem}\label{thm:hensel-smooth}
Let $k$ be a number field, $v$ a place of $k$, $\O_v$ the ring of integers in $k_v$ and $\F_v$ the residue field at $v$.  Let $X$ be a scheme over $\O_v$.  If $P \in X(\F_v)$ is a smooth point of $X \times_{\O_v} \F_v$, then $P$ lifts to a point of $X(\O_v)$.
\end{theorem}

Testing solubility of a variety over $\R$ is also a finite process: see~\cite[Theorem~13.13]{BPR:ARAG}.

To show that testing solubility at \emph{all} places of $\Q$ (that is, in $\R$ and in each $\Q_p$) is a finite procedure requires a little more work.  We would like to show that, for all primes $p$ above an explicitly determined bound, $X(\Q_p)$ is automatically non-empty.  
The main ingredient is the famous theorem of Lang and Weil on the number of points of a variety over a finite field.

\begin{theorem}[Lang--Weil~\cite{LW:AJM-1954}]\label{thm:lw}
There exists a constant $A(n,d,r)$ such that, for any finite field $\F$ of order $q$ and any geometrically irreducible variety $V \subset \P^n_\F$ of dimension $r$ and degree $d$, we have
\[
| \# V(\F) - q^r | < (d-1)(d-2) q^{r-\frac{1}{2}} + A(n,d,r) q^{r-1} .
\]
\end{theorem}

In particular, the theorem shows that there is a bound $B$, depending only on $n$, $r$ and $d$, such that $q>B$ implies $V(\F) \neq \emptyset$.  This bound can be explicitly determined (though doing so in general may be difficult).  If $V$ is smooth, then the Weil conjectures give a more detailed means to estimate the number of points on $V$.

The hypothesis of geometric irreducibility is essential: take $p>2$, let $a$ be any non-square in $\F_p$ and consider the variety $x^2-ay^2 \subset \A^2_{\F_p}$.  This only has the single point $(0,0)$, no matter how large $p$ is.  However, it is not geometrically irreducible since the equation factorises over $\F_{p^2}$.

Returning to our variety $X$ over $\Q$, Theorem~\ref{thm:lw} shows that, for $p$ sufficiently large, $X$ has a smooth $\F_p$-point.  (If $X$ is smooth, then so is $X \times \F_p$ for $p$ outside an explicitly computable set, so this is clear; otherwise, applying the theorem to the singular subscheme of $X$ shows that the number of singular points grows less quickly than the number of smooth points, so again existence of a smooth point is assured for $p$ sufficiently large.)  By Theorem~\ref{thm:hensel-smooth}, $X(\Q_p)$ is non-empty for $p$ sufficiently large.
Since only finitely many primes are left, testing solubility in all $\Q_p$ is seen to be a finite procedure.

The results described above generalise without problem to all number fields.

\subsection{Adelic points}

For ease of notation, we now define the set of adelic points of a variety over a number field.
Let $k$ be a number field and write $\pl_k$ for the set of places of $k$.  
For each $v \in \pl_k$, let $k_v$ be the completion of $k$ at $v$.  If $v$ is finite, let $\O_v$ be the ring of integers in $k_v$; if $v$ is infinite, set $\O_v = k_v$.
The ring of \emph{ad\`eles} of $k$ is the ring
\[
\ad_k = \{ (P_v) \in \prod_{v \in \pl_k} \mid P_v \in \O_v \text{ for almost all $v$} \}.
\]
(Here ``almost all'' means ``all but finitely many''.)
The ring $\ad_k$ comes with the restricted product topology.

Let $X$ be a variety over $k$.  
Since $\ad_k$ is a ring containing $k$, we can consider the set $X(\ad_k)$, and it turns out to be equal to
\[
X(\ad_k) = \{ (P_v) \in \prod_{v \in \pl_k} X(k_v) \mid P_v \in X(\O_v) \text{ for almost all $v$} \}.
\]
(Here we don't actually need to fix a model over the ring of integers of $k$, since any two models coincide at all but finitely many places.)

Each $X(k_v)$ comes equipped with a topology:  locally, $X(k_v)$ can be embedded in affine space $\A^n(k_v) = k_v^n$, so acquires a topology from that on $k_v$, and one can show that this does not depend on the choice of embedding.  The topology on $X(\ad_k)$ is then the restricted product topology with respect to the subsets $X(\O_v)$.  Explicitly, a basis of open sets is given by sets $U$ of the following form: take a finite set $S \subset \pl_k$, and for each $v \in S$ take an open set $U_v \in X(k_v)$; then let $U = \prod_{v \in S} U_v \times \prod_{v \notin S} X(\O_v)$.

If $X$ is projective, then we have $X(\O_v) = X(k_v)$: this is either an application of the valuative criterion of properness, or else the observation that any point of $\P^n(k_v)$ can be scaled to have coordinates lying in $\O_v$.  Hence, when $X$ is projective, we have $X(\ad_k) = \prod_v X(k_v)$, and the topology coincides with the product topology.

\subsection{The Hasse principle}

Let $X$ be a projective variety over a number field $k$.  We have seen that existence of points over all completions of $k$, that is, the condition $X(\ad_k) \neq \emptyset$, is a necessary condition for the existence of a $k$-rational point on $X$, and indeed a condition that can be checked by a finite procedure.
Sometimes this condition is also a sufficient condition, as in the following famous theorem of Hasse and Minkowski.
See~\cite[Chapter~IV]{Serre:CA} for a treatment over $\Q$.

\begin{theorem}[Hasse--Minkowski]\label{thm:hm}
Let $X \subset \P^n_k$ be a smooth hypersurface defined by one equation of degree $2$.  Then we have
\[
X(\ad_k) \neq \emptyset \quad \Longleftrightarrow \quad X(k) \neq \emptyset.
\]
\end{theorem}
This theorem means that testing whether such an $X$ admits a $k$-rational point is a finite process.
(But note that the theorem does not help you actually find a rational point.)

More generally, we define the Hasse principle as follows.

\begin{definition}
Let $X$ be a projective variety over a number field $k$.  We say that $X$ satisfies the \emph{Hasse principle} if the implication
\[
X(\ad_k) \neq \emptyset \quad \Longrightarrow \quad X(k) \neq \emptyset
\]
holds.
\end{definition}
This is frequently a more useful concept in the context of a class of varieties rather than for a single variety.  For example, Theorem~\ref{thm:hm} states that quadrics of any dimension satisfy the Hasse principle.

Other varieties known to satisfy the Hasse principle include Severi--Brauer varieties (that is, $X$ such that $X \times_k \kbar \cong \P^n_\kbar$ for some $n$) and del Pezzo surfaces of degree at least $5$.

However, not all varieties do satisfy the Hasse principle.  The following counterexample to the Hasse principle was discovered by Lind~\cite{Lind:thesis} and Reichardt~\cite{Reichardt:JRAM-1942}.

\begin{example}[Lind, Reichardt]\label{eg:lr}
Let $C$ be the smooth projective curve over $\Q$ defined in weighted projective space $\P(1,2,1)$ by the equation
\begin{equation}\label{eq:rl}
2Y^2 = X^4 - 17 Z^4.
\end{equation}
Then $C(\ad_\Q)$ is non-empty, but $C(\Q)$ is empty.
\end{example}
\begin{proof}
That $C(\ad_\Q)$ is non-empty is left as an exercise.  We prove that $C$ has no rational points.

If there were a rational point on $C$, then without loss of generality we could write it as
$(x,y,z)$ with $x,y,z$ integers and $x,z$ coprime.  What primes 
may divide $y$?  If $q>2$ is prime and $q
\divs y$, then we have $x^4 \equiv 17 z^4 \pmod q$ and so $17$ is a square
modulo $q$.  By quadratic reciprocity, this means that $q$ is a square
modulo $17$.

Since $2$ and $-1$ are also squares modulo $17$, we deduce that $y$ is 
a product of squares modulo $17$ and thus $y$ is a square modulo $17$.  We
can therefore write $y \equiv y_0^2 \pmod{17}$.  Substituting
into~\eqref{eq:rl}, we get $2y_0^4 \equiv x^4 \pmod{17}$ and
hence that $2$ is a fourth power modulo $17$.  But this is not true,
and so there can be no rational solution.
\end{proof}

Most of the arguments in this proof are entirely local arguments: they
involve making deductions about $X(\Q_p)$ for various places $p$.  But
there is one step which is not local, and that is the use of quadratic
reciprocity.  The theorem of quadratic reciprocity gives a link
between behaviour at one prime and behaviour at another prime, and thus
shows that the possible locations of our hypothetical rational
solution in the various $X(\Q_p)$ are not independent of each other.

Further counterexamples to the Hasse principle were discovered during the middle of the 20th century, including:
\begin{itemize}
\item a smooth diagonal cubic curve (Selmer~\cite{Selmer:AM-1951}, 1951);
\item a smooth cubic surface (Swinnerton-Dyer~\cite{SD:M-1962}, 1962), disproving a conjecture of Mordell that such surfaces satisfied the Hasse principle;
\item a smooth diagonal cubic surface (Cassels--Guy~\cite{CG:M-1966}, 1966);
\item a smooth del Pezzo surface of degree $4$ (Birch--Swinnerton-Dyer~\cite{BSD:JRAM-1975}, 1975).
\end{itemize}
In each case, the proof that the variety has no rational points used some form of reciprocity law.

\section{The Brauer--Manin obstruction}

In 1970, Manin~\cite{Manin:GBG} put forward a construction that explains all counterexamples to the Hasse principle known at that time.  This is what we now call the Brauer--Manin obstruction to the Hasse principle.
(At the time, it was not known that the Brauer--Manin obstruction explained the counterexample of Cassels--Guy mentioned above, but this was later shown to be true by Colliot-Th\'el\`ene--Kanevsky--Sansuc~\cite{CTKS}.)
More recently, further counterexamples to the Hasse principle have been found that are not explained by the Brauer--Manin obstruction, first by Skorobogatov~\cite{Skorobogatov:IM-1999}.  Despite this, the Brauer--Manin obstruction remains a central area of interest in the study of rational points, because it is conjectured to be the only obstruction to the Hasse principle for certain large classes of varieties.
Specifically, Colliot-Th\'el\`ene has conjectured this to be the case for geometrically rationally connected varieties.
Moreover, on many varieties the Brauer--Manin obstruction gives an explicitly computable obstruction to the Hasse principle.  On varieties for which it is known to be the only obstruction, this therefore gives a finite procedure for testing the existence of a rational point.

\subsection{Brauer groups of fields}

To define the Brauer--Manin obstruction, we first need to discuss Brauer groups of fields and of varieties.
An excellent reference for this material is the book by Gille and Szamuely~\cite{GS:CSAGC}.
The Brauer group of a field classifies certain objects called central simple algebras.

\begin{definition}
Let $K$ be a field.  A \emph{central simple algebra} over $K$ is a finite-dimensional $K$-algebra $A$ that is central (its centre is $K$) and simple (it has no two-sided ideals apart from $A$ and $0$).
\end{definition}

\begin{remark}
Our algebras always contain a $1$.  The $K$-multiples of $1$ form a copy of $K$ inside the algebra, so asking that the centre be equal to $K$ does make sense.
\end{remark}

Such an object is defined by a finite amount of data.  Indeed, a finite-dimensional $K$-algebra is a finite-dimensional $K$-vector space $V$ together with a $K$-linear multiplication law; to define the multiplication law, it is enough to choose a basis and specify the product of each pair of basis vectors.

\begin{example}
If $K$ is any field and $n \ge 1$ is an integer, then the matrix algebra $\Mat_n(K)$ is a central simple algebra over $K$ of dimension $n^2$.
\end{example}

\begin{example}
Take $K=\R$.  Hamilton's quaternion algebra $\HH$ is defined to be the $4$-dimensional $\R$-algebra with basis $1,i,j,k$ and multiplication given by $i^2=j^2=k^2=-1$, $ij=-ji=k$.
The algebra $\HH$ is a central simple algebra over $\R$, and indeed a division algebra.
\end{example}

\begin{example}
Let $K$ be any field of characteristic $\neq 2$, and let $a,b$ be two non-zero elements of $K$.  We can define the (generalised) quaternion algebra $(a,b)_K$ to be the $4$-dimensional vector space over $K$ with basis $1,i,j,k$ and multiplication given by $i^2=a$, $j^2=b$ and $ij=-ji=k$.  This is a central simple algebra over $K$; it is either a division algebra, or isomorphic to $\Mat_2(K)$.
\end{example}

We note here several properties (not all obvious) of quaternion algebras over a field that we will make use of later.  Assume that $K$ has characteristic different from $2$ and let $a,b,c$ be non-zero elements of $K$.
\begin{itemize}
\item The algebras $(a,b)_K$ and $(b,a)_K$ are isomorphic.
\item The algebras $(a,b)_K$ and $(a,bc^2)_K$ are isomorphic.
\item More generally, if $c$ lies in the norm group $\Norm_{K(\sqrt{a})/K} K(\sqrt{a})^\times$, then $(a,b)_K$ and $(a,bc)$ are isomorphic.
\item The algebra $(a,b)_K$ is isomorphic to $\Mat_2(K)$ if and only if $b$ lies in the norm group $\Norm_{K(\sqrt{a})/K} K(\sqrt{a})^\times$.
\end{itemize}

If $A,B$ are two central simple algebras over a field $K$, then the tensor product $A \otimes_K B$ is also a central simple algebra over $K$.  
In order to classify central simple algebras over $K$, we first put an equivalence relation on them that makes the matrix algebras ``trivial''.

\begin{definition}
Two central simple algebras $A,B$ over a field $K$ are \emph{equivalent} if there exist integers $m,n>0$ such that the algebras $A \otimes_K \Mat_m(K)$ and $B \otimes_K \Mat_n(K)$ are isomorphic.
\end{definition}

If $A$ is an algebra over $K$, then the \emph{opposite algebra} $A^\op$ is defined to have the same underlying vector space as $A$, but with multiplication defined in the opposite order.  There is a homomorphism $A \otimes_K A^\op \to \End_K(A)$ defined by sending $x \otimes y$ to the endomorphism $a \mapsto xay$, and it turns out that $A$ is a central simple algebra if and only if this homomorphism is an isomorphism.

We are now in a position to define the Brauer group of a field.

\begin{definition}
Let $K$ be a field.  The \emph{Brauer group} of $K$, written $\Br K$, is the group of equivalence classes of central simple algebras over $K$, with the operation of tensor product over $K$.
\end{definition}
Since $\End_K(A)$ is a matrix algebra, the discussion above shows that the class of $A^\op$ is the inverse of the class of $A$.

If $K \to L$ is a homomorphism of fields and $A$ is a central simple algebra over $K$, then $A \otimes_K L$ is a central simple algebra over $L$; this gives a homomorphism $\Br K \to \Br L$, making $\Br$ into a covariant functor.

We state without proof several standard facts about the Brauer groups of various fields.
\begin{itemize}
\item If $K$ is separably closed, then $\Br K$ is trivial.
\item If $\F$ is a finite field, then $\Br \F$ is trivial.
\item $\Br \R$ has order $2$; the trivial class is $\{\Mat_n(\R) : n \in \N \}$ and the non-trivial class is $\{ \Mat_n(\HH) : n \in \N \}$.
\item The class of a quaternion algebra over $K$ is always of order $2$ in $\Br K$, but not every class of order $2$ contains a quaternion algebra.  (For example, the tensor product of two quaternion algebras is not necessarily equivalent to a quaternion algebra.)
\item If $K$ is a finite extension of $\Q_p$, then there is a canonical isomorphism $\inv \colon \Br K \to \Q/\Z$ called the Hasse invariant map.
In particular, for the quaternion algebra $(a,b)_K$, we have
\[
\inv(a,b)_K = \begin{cases}
0 & \text{if $b \in \Norm_{K(\sqrt{a})/K} K(\sqrt{a})^\times$}; \\
\frac{1}{2} & \text{otherwise}.
\end{cases}
\]
\item Let $k$ be a number field. For each finite place $v$ of $k$, let $\inv_v \colon \Br k_v \to \Q/\Z$ be the Hasse invariant map.  Moreover, for each infinite place $v$, define $\inv_v \colon \Br k_v \to \Q/\Z$ be the unique injective homomorphism (with image of order $2$ if $v$ is a real place, and trivial image if $v$ is a complex place).
Then the sequence
\begin{equation}\label{eq:abhn}
0 \to \Br k \to \bigoplus_{v \in \pl_k} \Br k_v \xrightarrow{\sum_v \inv_v} \Q/\Z \to 0
\end{equation}
is a short exact sequence.  This is the celebrated theorem of Albert--Brauer--Hasse--Noether.  Implicit in here is the statement that, for $A \in \Br k$, the image of $A$ in $\Br k_v$ is zero for all but finitely many $v$.  Note that it really matters that we use the canonical map $\inv_v$ and not just any isomorphism $\Br k_v \to \Q/\Z$.
\end{itemize}

For any field $K$ with separable closure $K^s$, it can be shown that $\Br K$ is canonically isomorphic to the Galois cohomology group $\H^2(K, (K^s)^\times)$.

The exact sequence~\eqref{eq:abhn}, when restricted to classes of quaternion algebras, is equivalent to quadratic reciprocity.  In general, \eqref{eq:abhn} can be seen as a generalised reciprocity law, and it is this on which the Brauer--Manin obstruction is based.

\subsection{The Brauer--Manin obstruction}

In the next section, we will define the Brauer group of a variety.  For now, we use a vague notion.
Let $X$ be a variety over a number field $k$.  An element of the Brauer group $\Br X$ will be a ``family of Brauer group elements'' parametrised by $X$.  More specifically, given $\al \in \Br X$ and a point $P \in X(L)$ for some field $L \supset K$, one can ``evaluate'' $\al$ at $P$ to obtain $\al(P) \in \Br L$.  This property is enough to define the obstruction.

Fix an element $\al \in \Br X$.  Given an adelic point $(P_v) \in X(\ad_k)$, we can evaluate $\al$ at each $P_v$ separately and obtain an element of $\prod_v \Br k_v$.  This turns out to lie in $\bigoplus_v \Br k_v$, and so we can use the exact sequence~\eqref{eq:abhn} to test whether it might come from $\Br k$.  If not, then $(P_v)$ certain does not come from a point $P \in X(k)$.

Putting this idea into action gives a commutative diagram as follows.
\begin{center}
\begin{tikzcd}
X(k) \arrow[r] \arrow[d, "\al"'] & X(\ad_k) \arrow[d, "\al"'] \\
\Br k \arrow[r] & \bigoplus_v \Br k_v \arrow[r, "\sum_v \inv_v"] & \Q/\Z.
\end{tikzcd}
\end{center}

We define a subset $X(\ad_k)^\al$ of the adelic points of $X$, consisting of those which pass this test.  Specifically, define
\[
X(\ad_k)^\al = \{ (P_v) \in X(\ad_k) \mid \sum_v \inv_v \al(P_v) = 0 \}
\]
and
\begin{align*}
X(\ad_k)^{\Br} &= \bigcap_{\al \in \Br X} X(\ad_k)^\al \\
&= \{ (P_v) \in X(\ad_k) \mid \sum_v \inv_v \al(P_v) = 0 \text{ for all } \al \in \Br X \}.
\end{align*}

The diagram above shows that the diagonal image of $X(k)$ in $X(\ad_k)$ is contained in $X(\ad_k)^{\Br}$.  It follows that, if $X(\ad_k)^{\Br}$ is empty, then there are no rational points on $X$.

\begin{definition}
If $X(\ad_k)$ is non-empty but $X(\ad_k)^{\Br}$ is empty, we say that there is a \emph{Brauer--Manin obstruction to the Hasse principle} on $X$.
\end{definition}

It can be shown that the map $X(k_v) \to \Br k_v$ given by an element of $\Br X$ is continuous; so $X(\ad_k)^{\Br}$ is a closed subset of $X(\ad_k)$.

\subsection{Weak approximation}

We conclude this section by mentioning the related idea of weak approximation.  Using the diagonal map, we can think of $X(k)$ as a subset of $X(\ad_k)$.

\begin{definition}
Let $X$ be a variety over the number field $k$.  Then $X$ satisfies \emph{weak approximation} if $X(k)$ is dense in $\prod_v X(k_v)$.
\end{definition}

It is helpful to unpack this definition using the definition of the product topology.
A basic open set $U \subset \prod_v X(k_v)$ is given by a finite set $S \subset \pl_k$ and open sets $U_v \subset X(k_v)$ for $v \in S$; then $U$ is the product $\prod_{v \in S} U_v \times \prod_{v \notin S} X(k_v)$.
That $X(k)$ be dense in $\prod_v X(k_v)$ is the same as that every such basic open set contain a point of $X(k)$.
An open set of $X(k_v)$, for $v$ a finite place, is defined by congruence conditions; so saying that $X$ satisfies weak approximation is a little like a ``Chinese remainder theorem'' for rational points of $X$:  given congruence conditions at finitely many places, there exists a point of $X(k)$ satisfying those conditions.

The Brauer--Manin obstruction can also obstruct weak approximation on a projective variety $X$.  Given that $X(\ad_k)^{\Br}$ is closed in $X(\ad_k)$, it is impossible for $X$ to satisfy weak approximation as soon as $X(\ad_k)^{\Br}$ is not the whole of $X(\ad_k)$.

\begin{definition}
If $X(\ad_k)$ is non-empty and $X(\ad_k)^{\Br}$ is not the whole of $X(\ad_k)$, then we say there is a \emph{Brauer--Manin obstruction to weak approximation} on $X$.
\end{definition}

\section{Understanding Brauer groups}

\subsection{The Brauer group of a variety}

Let $X$ be a smooth, projective, geometrically irreducible variety over a number field $k$.
Given a central simple algebra $\al$ over the function field $k(X)$, we can try to evaluate $\al$ at points of $X$ in the obvious way.

\begin{example}
Take $X = \A^1_\Q$, with function field $\Q(t)$, and $\al = (-1,t)_{\Q(t)}$.  For any $a \neq 0$, substituting $t=a$ gives a quaternion algebra $(-1,a)_\Q$ over $\Q$.  For $a=0$, this does not work.
\end{example}

\begin{example}
Let $X$ be the elliptic curve over $\Q$ defined by the Weierstrass equation
\[
y^2=(x-e_1)(x-e_2)(x-e_3)
\]
and let $\al$ be the quaternion algebra $(3,x-e_1)_{\Q(X)}$.  Evaluating at any point of $X(\Q)$ apart from $(e_1,0)$ gives a well-defined quaternion algebra over $\Q$, but it looks as if we cannot evaluate $\al$ at the point $(e_1,0)$.
However, dividing $x-e_1$ by $y^2$ shows that $\al$ is isomorphic to the algebra $(3,((x-e_2)(x-e_3))^{-1})_{\Q(X)}$, and here we can indeed substitute $(e_1,0)$ to get an algebra over $\Q$.  So it seems that, just as with rational functions, it is sometimes possible to extend the domain of definition by writing the algebra in a different way.
\end{example}

In the previous example we replaced our algebra by an isomorphic one in order to evaluate it at a particular point.
In general we will be interested in evaluating not an algebra but a class in $\Br k(X)$ at a point of $X$, so we can replace our algebra by an equivalent one if necessary.

As it stands, our idea of ``evaluating'' at a point $P \in X$ is unsatisfying because it seems to rely on an evaluation homomorphism $k(X) \to k(P)$, which of course does not exist.  
The way to understand this correctly for rational functions is to use the local ring at $P$.  The ring $\O_{X,P} \subset k(X)$ consists of those functions which are defined at $P$, and there is an evaluation homomorphism $\O_{X,P} \to k(P)$.
To apply the same idea to evaluating Brauer group elements, we must define the Brauer group of the local ring $\O_{X,P}$.

\begin{definition}
Let $R$ be a local ring with maximal ideal $\m$.  An \emph{Azumaya algebra} over $R$ is a free $R$-algebra $A$ such that $A/\m A$ is a central simple algebra over the field $R/\m$.  Two Azumaya algebras $A,B$ over $R$ are \emph{equivalent} if there exist $m,n>0$ such that $A \otimes_R \Mat_m(R)$ and $B \otimes_R \Mat_n(R)$ are isomorphic.  The \emph{Brauer group} of $R$ consists of the equivalence classes of Azumaya algebras over $R$, under the operation of tensor product over $R$.
\end{definition}

If $R$ has field of fractions $K$, then there is a natural map $\Br R \to \Br K$ given by sending an Azumaya algebra $A$ to the central simple algebra $A \times_R K$.
Auslander and Goldman proved that, if $R$ is a regular local ring, then $\Br R \to \Br K$ is an injection~\cite[Theorem~7.2]{AG:TAMS-1960}.

Now let $P$ be a point of a smooth variety $X$ over a field $k$.  The ring $\O_{X,P}$ is regular, so we can regard $\Br \O_{X,P}$ as a subgroup of $\Br k(X)$.  If a class $\al \in \Br k(X)$ lies in $\Br \O_{X,P}$, then we can use the natural map $\Br \O_{X,P} \to \Br k(P)$ to evaluate $\al$ at $P$.  (This evaluation does indeed coincide with the na\"ive notion used above.)

We would like $\Br X$ to consist of elements of $\Br k(X)$ that can be evaluated at all points of $X$, so it makes sense to define $\Br X$ as follows.

\begin{definition}
Let $X$ be a smooth variety over a field $k$.  The \emph{Brauer group} of $X$ is the group
\[
\Br X = \bigcap_{P \in X} \Br \O_{X,P},
\]
the intersection being taken inside $\Br k(X)$.
\end{definition}

Note that the intersection here is over all scheme-theoretic points of $X$, though it is straightforward to show that one may equivalently take all closed points.  (If a point $P$ is a specialisation of another point $Q$, then we have $\O_{X,P} \subset \O_{X,Q} \subset k(X)$ and therefore $\Br \O_{X,P} \subset \Br \O_{X,Q} \subset \Br k(X)$.)

This definition of the Brauer group has the merit of being fairly elementary.  Another common definition of the Brauer group is that it is the \'etale cohomology group $\H^2_\et(X,\Gm)$; this is useful if one wants to use the powerful tools of \'etale cohomology to study the Brauer group.
In the case of a smooth, projective variety, the two definitions coincide: this follows from a difficult result of Gabber, of which an alternative proof was given by de Jong~\cite{DeJong:Gabber}.

\subsection{An example of the Brauer--Manin obstruction}

Now that we have defined the Brauer group of a variety, we are in a position to use the Brauer--Manin obstruction to explain the counterexample to the Hasse principle given in Example~\ref{eg:lr}.

Let $C$ be the curve $2Y^2 = X^4-17Z^4$ as before.  We will use the quaternion algebra $\al=(17,Y/X^2) \in \Br \Q(C)$, which actually lies in $\Br C$.  Note that $\al$ is isomorphic to the algebra $(17,Y/Z^2)$ so we may use whichever representation is more convenient.

To evaluate the Brauer--Manin obstruction, we must compute the value of the invariant $\inv_v \al(P_v)$ for every place $v$ and every point $P_v \in X(\Q_v)$.  We look at the places individually.

If $17$ is a square in $\Q_v$, then we have $\inv_v \al(P_v)=0$ for all $P_v \in X(k_v)$.  As long as either $Y/X^2$ or $Y/Z^2$ is invertible at $P_v$, we can see this by first substituting the coordinates of $P$ into one of the expressions for $\al$ and seeing that the resulting quaternion algebra is trivial over $k_v$.
However, we can also see this more efficiently by noting that the evaluation map $\Br C \to \Br k_v$ factors through $\Br(C \times_k k_v)$, and that $\al$ is zero in this group since $17$ is a square in the function field $k_v(C)$.
Note that this case includes the places $2$ and $\infty$.

If $v=p$ is an odd prime other than $17$, and $17$ is not a square in $\Q_p$, then we proceed as follows.  Any point $P_v$ of $X(\Q_p)$ can be written as $(x,y,z)$ with $x,y,z \in \Z_p$ not all divisible by $p$.  If $p$ were to divide $y$, then looking at the equation modulo $p$ shows that $17$ would be a square modulo $p$, which is not true; so $p$ does not divide $y$.  It follows that at least one of $y/x^2$ and $y/z^2$ is a unit in $\Z_p$.
The extension $\Q_p(\sqrt{17})/\Q_p$ is unramified, so every unit in $\Z_p$ is a norm from this extension, giving $\inv_v \al(P_v)=0$.

Finally, take $v=17$ and let $P=(x,y,z)$ be a point of $C(\Q_{17})$.  Again we may assume that $x,y,z$ lie in $\Z_{17}$.  If, say, $x$ is divisible by $17$, then looking modulo $17$ shows $17 \divs y$, but then $17^2 \divs 17z^4$ and so $17 \divs z$, a contradiction.  So $x$ is not divisible by $17$, and similarly neither is $y$.  Now looking modulo $17$, and using the fact that $2$ is not a fourth power in $\F_{17}$, shows that $y$ is not a square modulo $17$.  The elements of $\Z_{17}^\times$ that are norms from $\Q_{17}(\sqrt{17})$ are precisely those which are square modulo $17$, so we deduce $\inv_{17} \al(P) = \frac{1}{2}$.

To summarise, we have found $\inv_{17} \al(P)=\frac{1}{2}$ for all points $P \in C(\Q_{17})$, and $\inv_v \al(P_v)=0$ for all $v \neq 17$ and $P_v \in C(\Q_v)$.
Therefore, for all $(P_v) \in C(\ad_\Q)$, we have $\sum_v \inv_v \al(P_v) = \frac{1}{2}$ and therefore $C(\ad_\Q)^{\Br} = \emptyset$, showing that there is a Brauer--Manin obstruction to the Hasse principle on $C$.

\begin{remark}
The ingredients for computing the Brauer--Manin obstruction on $C$ were exactly the same as the ingredients in our earlier proof that $C$ has no rational point.
\end{remark}

\subsection{Brauer groups of some varieties}

Here we collect some facts about Brauer groups of varieties.

If $k$ is algebraically closed, then we have $\Br(\P^n_k)=0$ for all $n$, and also $\Br(C)=0$ for every curve $C$ over $k$.  (The latter follows from Tsen's theorem that $\Br k(C)$ is a $C_1$ field and so has trivial Brauer group.)
If $X$ is a smooth proper variety over an algebraically closed field $k$ of characteristic zero, then $\Br X$ contains $(\Q/\Z)^{b_2-\rho}$, where $b_2$ is the second Betti number of $X$ and $\rho$ is the rank of the N\'eron--Severi group of $X$.

Over any field $k$ of characteristic zero, the Brauer group is a birational invariant of proper smooth varieties: if $X$ and $Y$ are birational proper smooth varieties, then $\Br X$ and $\Br Y$ are isomorphic.

Let $X$ be a smooth, proper, geometrically irreducible variety over a field $k$.
Let $\Xbar$ denote the base change of $X$ to a separable closure of $k$.
We define some subgroups of $\Br X$.
\begin{itemize}
\item The subgroup $\Br_0 X = \im(\Br k \to \Br X)$ consists of the ``constant'' elements.  (Beware: the homomorphism $\Br k \to \Br X$ is not necessarily injective, though it is if $X$ has a $k$-point.)
\item The subgroup $\Br_1 X = \ker(\Br X \to \Br \Xbar)$ consists of the ``algebraic'' elements, that is, those that become trivial after a finite extension of the base field $k$.  (This clearly contains $\Br_0 X$.)
\item The quotient $\Br X / \Br_1 X \cong \im(\Br X \to \Br \Xbar)$ is called the ``transcendental'' Brauer group of $X$.
\end{itemize}

Putting these facts and definitions together allows us to make some new statements, such as: if $X$ is a geometrically rational variety, then $\Br \Xbar$ is trivial and so $\Br_1 X$ is the whole of $\Br X$.

Let $X$ be smooth, proper and geometrically irreducible over $k$.
A very useful calculation, coming from the Hochschild--Serre spectral sequence in \'etale cohomology, is the existence of the following exact sequence of Galois cohomology groups:
\[
0 \to \Pic X \to \H^0(k, \Pic \Xbar) \to \Br k \to \Br_1 X \to \H^1(k, \Pic \Xbar) \to \H^3(k,\kbar^\times).
\]
If $k$ is a number field, then we have $\H^3(k,\kbar^\times)=0$ and so an extremely useful isomorphism
\begin{equation}\label{eq:h1pic}
\Br_1 X / \Br_0 X \cong \H^1(k, \Pic\Xbar).
\end{equation}
Supposing that $\Pic\Xbar$ is finitely generated and torsion-free (for example, $X$ is a rational surface or a K3 surface), the group $\H^1(k,\Pic\Xbar)$ is finite and can easily be computed explicitly.
Inverting the isomorphism to obtain explicit elements of $\Br_1 X$ is more challenging, but in practice can be done in many examples.

Note that elements of $\Br_0 X$ give no Brauer--Manin obstruction, and the obstruction given by $\al \in \Br X$ depends only on the class of $\al$ modulo $\Br_0 X$, so in practice we are usually interested in the quotient $\Br X / \Br_0 X$.

\subsection{Residue maps}

Let $X$ be a smooth variety over a field $k$.
Given an element of $\Br k(X)$, how can one tell whether it lies in $\Br X$?  The corresponding question for rational functions has a pleasing answer:  a rational function is regular if and only if it has trivial valuation along every prime divisor of $X$.
In this section we describe the corresponding result for Brauer groups.  The first aim is to define a residue map for each prime divisor.

Let $R$ be a discrete valuation ring with maximal ideal $\m$, field of fractions $K$ and perfect residue field $F$.
\begin{theorem}\label{thm:res}
There is a homomorphism $\partial \colon \Br K \to \H^1(F,\Q/\Z)$ whose kernel is precisely $\Br R$.
\end{theorem}
The homomorphism $\partial$ is called the \emph{residue map} associated to the DVR $R$.

We briefly describe the construction of the residue map.  Firstly, assume that $K$ is complete.
It can be shown that every element of $\Br K$ is split by some finite unramified extension of $K$.
Letting $K^\nr$ be the maximal unramified extension of $K$, this shows that the natural map $\Br K \to \Br K^\nr$ is the zero map; the inflation-restriction sequence in Galois cohomology (together with Hilbert's Theorem 90) gives an exact sequence
\[
0 \to \H^2(K^\nr/K, (K^\nr)^\times) \to \Br K \to \Br K^\nr,
\]
showing that $\H^2(K^\nr/K,(K^\nr)^\times) \to \Br K$ is an isomorphism.

The valuation map on $K$ extends to a homomorphism $(K^\nr)^\times \to \Z$.  This induces a homomorphism
\[
\H^2(K^\nr/K,(K^\nr)^\times) \to \H^2(K^\nr/K, \Z) 
\cong \H^2(F,\Z) \cong \H^1(F,\Q/\Z).
\]
The first isomorphism comes from the isomorphism $\Gal(K^\nr/K) \cong \Gal(\bar{F}/F)$, because $K$ is complete, and the second from the boundary map in the long exact sequence in cohomology coming from $0 \to \Z \to \Q \to \Q/\Z \to 0$.
This composition is the residue map.

If $K$ is not complete, let $\hat{K}$ be its completion, and define the residue map via the natural map $\Br K \to \Br \hat{K}$.

\begin{remark}\label{inv}
When $K$ is a finite extension of $\Q_p$, there is a canonical isomorphism $\H^1(\F_q,\Q/\Z) \to \Q/\Z$ given by evaluation at the Frobenius element.  Composing the residue map with this isomorphism gives the Hasse invariant map $\inv$ (and indeed is often given as the definition of the invariant map).
\end{remark}

\begin{remark}
The condition that the residue field be perfect is essential.  For example, take $K=\Q_2(t)$ with the $2$-adic valuation, with residue field $F=\F_2(t)$.  The quaternion algebra $(2,t)$ has no splitting field that is unramified over $K$; of the two obvious ones, $K(\sqrt{2})$ is wildly ramified, and $K(\sqrt{t})$ is even worse, inducing an inseparable extension of $F$.
\end{remark}

Returning to the case of a smooth variety $X$ over a field $k$ of characteristic zero, we can apply Theorem~\ref{thm:res} to the local ring at every prime divisor and so obtain, for each prime divisor $D \subset X$, a residue map $\partial_D \colon \Br k(X) \to \H^1(k(D),\Q/\Z)$.  One reason that these residue maps are so useful is the purity theorem.

For $\al \in \Br k(X)$, define the ramification locus of $\al$ to be the set 
\[
\{ P \in X \mid \al \notin \Br \O_{X,P} \}.
\]
This is a proper closed subset of $X$:  one way to see this is to use the isomorphism $\al \otimes_{k(X)} \al^\op \to \End_{k(X)} \al$, and show that set of points where it fails to specialise to an isomorphism is given by the vanishing of a determinant.

\begin{theorem}[Purity theorem, Grothendieck, first form]
Let $X$ be a smooth variety over a field $k$ of characteristic zero.
For $\al \in k(X)$, the ramification locus of $\al$ is either empty or of pure codimension $1$ in $X$.
\end{theorem}

This allows us to give the following alternative description of the Brauer group of $X$:
\[
\Br X = \bigcap_D \Br \O_{X,D},
\]
where the intersection is over all prime divisors $D$ on $X$.

Using Theorem~\ref{thm:res}, we can give a second form of the purity theorem:

\begin{theorem}[Purity theorem, second form]\label{thm:purity}
Let $X$ be a smooth variety over a field $k$ of characteristic zero.
The following sequence is exact:
\[
0 \to \Br X \to \Br k(X) \xrightarrow{\oplus_D \partial_D} \bigoplus_D \H^1(k(D),\Q/\Z),
\]
where $D$ runs over all prime divisors on $X$.
\end{theorem}
Note that the map $\oplus_D \partial_D$ does indeed land in the direct sum: for each element $\al \in \Br k(X)$, we can have $\partial_D(\al) \neq 0$ only for the finitely many divisors $D$ lying in the ramification locus of $\al$.

\begin{remark}
The result also holds in positive characteristic, as long as one restricts attention to the prime-to-$p$ torsion in the Brauer group.
\end{remark}

\begin{remark}
If $D$ is itself smooth and $\al$ is unramified away from $D$, that is, $\al$ lies in $\Br(X \setminus D)$, then one can show that the residue $\partial_D(\al)$
lies in the subgroup $\H^1(D,\Q/\Z) \subset \H^1(k(D),\Q/\Z)$.
\end{remark}

\section{Local solubility in families}\label{sec:ls}

In this section we turn our attention from individual varieties to families of varieties.
The aim is to understand results such as those of Loughran and Smeets~\cite{LS:few} which relate the number of everywhere locally soluble varieties in a family to the geometry of the family.

\subsection{Introduction}

Let $k$ be a number field.  By a ``family of varieties'' over $k$ we mean a smooth, projective, irreducible variety $X$ over $k$ equipped with a dominant morphism $\pi \colon X \to \P^n_k$, such that the generic fibre of $\pi$ is geometrically irreducible.  (This condition ensures that the other fibres of $\pi$ will also be geometrically irreducible outside a proper closed subset of the base.)  We think of $X$ as a family of varieties in the sense that every $P \in \P^n(k)$ gives a $k$-variety $X_P = \pi^{-1}(P)$.

\begin{example}\label{eg:conics}
Take the ambient space $\P^2_k \times \P^2_k$ with coordinates $(x_0,x_1,x_2)$ on the first factor and $(a_0,a_1,a_2)$ on the second factor.  Let the subvariety $X$ be defined by the equation $a_0 x_1^2 + a_1 x_1^2 + a_2 x_2^2 = 0$.  If we take $\pi$ to be the projection to the second factor, restricted to $X$, then we obtain a family of varieties in which every $P = (a_0:a_1:a_2) \in \P^2(k)$ corresponds to the plane conic defined by the equation above.
\end{example}

The main question in this section will be: given such a family of varieties, how many of the varieties are everywhere locally soluble, that is, for how many of the points $P \in \P^n(k)$ do we have $X_P(\ad_k) \neq \emptyset$?
As it stands, this question does not make much sense: there are infinitely many points $P \in \P^n(k)$, and in most cases infinitely many of the corresponding $X_P$ will be everywhere locally soluble.  To turn the question into a more precise one, we use the standard height on $\P^n(k)$.

In the case $k=\Q$, we can define the height of a point $P \in \P^n(\Q)$ as follows: write $P = (a_0 : \dotsb : a_n)$ with $a_0, \dotsc, a_n$ coprime integers, and define the height to be $H(P) = \max\{ |a_0|, \dotsc, |a_n|\}$.
In the case of a more general number field, we may not be able to write a point using coprime integers, so we use the more general definition
\[
H((a_0:\dotsb:a_n)) = \prod_{v \in \pl_k} \max\{|a_0|_v, \dotsc, |a_n|_v \}.
\]
It is an easy but worthwhile exercise to check that this gives the same answer as above in the case $k=\Q$.

Given a family $\pi \colon X \to \P^n$ as above, we now define for $B>0$
\[
N_\mathrm{loc}(B) = \# \{ P \in \P^n(k) \mid H(P) \le B \text{ and } X_P(\ad_k) \neq \emptyset \}.
\]
The aim is to compare this to the total number of points of height at most $B$,
\[
N_\mathrm{tot}(B) = \# \{ P \in P^n(k) \mid H(P) \le B \}.
\]
A precise asymptotic for $N_\mathrm{tot}(B)$ is given by the following theorem of Schanuel.
\begin{theorem}[Schanuel]
Let $k$ be a number field.  Let $h$ be the class number of $k$, $\zeta_k$ the zeta function of $k$, $r_1$ and $r_2$ the numbers of real and complex places respectively, $r=r_1+r_2-1$ the rank of the unit group, $R$ the regulator, $\Delta$ the discriminant, $w$ the number of roots of unity in $k$.
Then we have
\[
N_\mathrm{tot}(B) \sim \frac{h}{\zeta_k(n+1)}
\Big(\frac{2^{r_1}(2\pi)^{r_2}}{\sqrt{|\Delta|}}\Big)^{n+1}
\frac{R}{w} (n+1)^r B^{(n+1)}
\]
\end{theorem}
(See Serre~\cite[Section~2.5]{Serre:LMWT}; the version stated there has $B^{(n+1)[k:\Q]}$ instead of $B^{n+1}$, but that is because it uses a different normalisation of the absolute values on $k$.)
The new question is: how does the ratio $N_\mathrm{loc}(B) / N_\mathrm{tot}(B)$ behave as $B \to \infty$?  Whenever we talk about a ``proportion'' of points in $\P^n(k)$ in future, it will be in this sense.

Let us return to Example~\ref{eg:conics} and try to answer the simpler question of solubility at just one place, beginning with solubility in $\R$.  
For $P =(a_0:a_1:a_2)$, it is easy to see that $X_P(\R)$ is non-empty if and only if $a_0,a_1,a_2$ do not all have the same sign.
This clearly happens for $3/4$ of all triples $a_0,a_1,a_2$ up to any given bound, but let us try to explain this in terms of the measure of some subset of $\P^2(\R)$.
Temporarily, let $a_0,a_1,a_2$ take values in $\R$.  If we scale $(a_0,a_1,a_2)$ to lie on the unit sphere (which is compact), then the two octants making up the points for which the corresponding conic is soluble make up precisely $3/4$ of the surface area of the sphere.
For any $B>0$, the points of $\P^2(\Q)$ of height $\le B$ are evenly distributed on the unit sphere in the sense that equal numbers lie in each octant, so precisely $3/4$ correspond to soluble fibres.
Of course we should really be working in $\P^2(\R)$, which is the quotient of the unit sphere by the antipodal map.

For a more general family of varieties, it is not reasonable to expect the proportion of points of height at most $B$ having soluble fibres to be independent of $B$; but the fact that $\P^n(\Q)$ is equidistributed in $\P^n(\R)$ as $B \to \infty$ means that the proportion of soluble fibres does tend to the ratio calculated using the measure of the relevant part of $\P^n(\R)$.

If we now turn to looking at solubility at an odd prime $p$, the situation is more complicated but the same basic idea applies.  Let $(a_0:a_1:a_2)$ be a point of $\P^2(\Q_p)$, scaled so that $a_0,a_1,a_2$ lie in $\Z_p$ and are not all divisible by $p$.
What are the conditions on $(a_0,a_1,a_2)$ for $a_0x_0^2 + a_1x_2^2 + a_2x_2^2$ to be soluble in $\Q_p$?

If none of the $a_i$ are divisible by $p$, then the equation is soluble: the Chevalley--Warning theorem~\cite[Chapter~I, Section~2, Theorem~3]{Serre:CA} shows that there is a solution over $\F_p$, which lifts by Hensel's Lemma.

If, say, $v_p(a_0)=1$ but $p$ does not divide $a_1a_2$, then the solubility depends on whether $a_1x_1^2 + a_2x_2^2$ has a non-zero solution modulo $p$, which happens if and only if $-a_1a_2$ is a square in $\F_p$.

The other cases look more complex but can be reduced to these by transformations: for example, if $p^2$ divides $a_0$ then substituting $x_0 \mapsto p x_0$ shows that this equation is soluble if and only if that corresponding to $(a_0/p^2:a_1:a_2)$ is soluble.
By keeping track of the effect of these transformations on the measure on $\P^2(\Q_p)$ we can produce an explicit expression for the proportion of $\P^2(\Q_p)$ corresponding to fibres soluble in $\Q_p$.
Using the fact that points of $\P^2(\Q)$, ordered by height, are equidistributed in $\P^2(\Q_p)$ shows that this is also the proportion of $\P^2(\Q)$ corresponding to fibres locally soluble in $\Q_p$.

So the proportion of varieties in our family that are locally soluble at one given place can be computed as the measure of some set.  What about the proportion that are locally soluble at all places?  This particular example was studied by Serre~\cite[Exemple~3]{Serre:CRAS-1990}.

\begin{theorem}[Serre]\label{thm:serre}
For the family of Example~\ref{eg:conics}, we have
\[
N_\mathrm{loc}(B) \ll \frac{B^3}{(\log B)^{3/2}}.
\]
In particular, $N_\mathrm{loc}(B) / N_\mathrm{tot}(B)$ tends to $0$ as $B \to \infty$.
\end{theorem}

On the other hand, there are families of varieties such that the ratio is positive.

\begin{theorem}[Bright, Browning, Loughran~\cite{BBL:wa}]\label{thm:cubics}
Let $X \subset \P^3 \times \P^3$ be the family of varieties over $\Q$ defined by
\[
a_0 x_0^3 + a_1 x_1^3 + a_2 x_2^3 + a_3 x_3^3 = 0.
\]
Then the ratio $N_\mathrm{loc}(B) / N_\mathrm{tot}(B)$ tends to approximately $0.86$ as $B \to \infty$.
Moreover, this limit is a product $\prod_{v \in \pl_k} c_v$, where $c_v$ is the proportion of varieties in the family that are soluble over $\Q_v$.
\end{theorem}

\begin{theorem}[Poonen, Voloch~\cite{PV:PM-2004}]
Fix positive integers $n,d \ge 2$ satisfying $(n,d) \neq (2,2)$.
For the family of all hypersurfaces of degree $d$ in $\P^n$, the ratio $N_\mathrm{loc}(B) / N_\mathrm{tot}(B)$ tends to a limit $c>0$ as $B \to \infty$.
Moreover, we have $c = \prod_{v \in \pl_k} c_v$, where $c_v$ is the proportion of varieties in the family that are soluble over $\Q_v$.
\end{theorem}

What is the qualitative difference between the two types of behaviour?

\subsection{Bad reduction and local solubility}

At the root of the relationship between geometry and local solubility of varieties is the result of Lang--Weil (Theorem~\ref{thm:lw}) which implies, approximately, that a geometrically irreducible variety over a sufficiently large finite field must have a point.

In the case of a scheme with several components, we of course only need one component to be geometrically irreducible in order to apply the Lang--Weil theorem.   Skorobogatov made the following definition.

\begin{definition}
A scheme of finite type over a field is \emph{split} if it contains a geometrically integral open subscheme.
\end{definition}

Equivalently, a scheme of finite type over a field is split if and only if it has an irreducible component that is geometrically irreducible and has multiplicity $1$.
This concept turns out to be very useful in studying the local solubility of families of varieties.

Let us look at the two contrasting examples from above, the diagonal cubic surfaces and the plane conics.

\begin{example}
Let $X \subset \P^3 \times\P^3$ be the family of varieties over $\Q$ defined by
\[
a_0 x_0^3 + a_1 x_1^3 + a_2 x_2^3 + a_3 x_3^3 = 0.
\]

Let $D \subset \P^3$ be the divisor defined by $a_0=0$, and let $d$ be the generic point of $D$.  The fibre $X_d$ is the variety in $\P^3$ over the field $\Q(a_1/a_3, a_2/a_3)$ defined by the equation
\[
(a_1/a_3) x_1^3 + (a_2/a_3) x_2^2 + x_3^3 = 0.
\]
Since the variable $x_0$ does not appear in the equation, this variety is a cone over a smooth cubic curve, with the vertex of the cone at $(1:0:0:0)$.  In particular, it is geometrically irreducible and therefore split.

Suppose that $P=(a_0:a_1:a_2:a_3) \in \P^3(\Q)$ satisfies $a_0a_1a_2a_3 \neq 0$, with $a_0,a_1a_2,a_3$ coprime integers, and that $p$ is an odd prime dividing $a_0$.  Then the reduction of $P$ modulo $p$ lands in the divisor $D$.  Assuming this happens in a sufficiently generic way (that is, $p^2 \ndivs a_0$ and $p \ndivs a_1a_2a_3$), then the reduction of $X_P$ modulo $p$ looks exactly like the fibre $X_d$, that is, a cone over a smooth cubic curve.
If $p$ is sufficiently large, then Lang--Weil implies that there is a smooth point on the reduction of $X_P$ modulo $p$, and therefore $X_P(\Q_p) \neq \emptyset$.
\end{example}

\begin{example}
Let $X \subset \P^2 \times \P^2$ be the family of varieties from Example~\ref{eg:conics}, defined by
\[
a_0x_0^2 + a_1x_1^2 + a_2x_2^2 = 0.
\]

As before, let $D \subset \P^2$ be the divisor defined by $a_0=0$ and let $d$ be the generic point of $D$.
The fibre $X_d$ is the variety in $\P^2$ over the field $K=\Q(a_1/a_2)$ defined by
\[
(a_1/a_2) x_2^2 + x_3^2 = 0.
\]
This is not split: geometrically, it is the union of two lines, conjugate over the extension field $K(\sqrt{-a_1/a_2})$.

Now suppose that $P = (a_0:a_1:a_2) \in \P^2_\Q$ satisfies $a_0a_1a_2 \neq 0$, with $a_0,a_1,a_2$ coprime integers, and that $p$ is an odd prime dividing $a_0$.  As before, the reduction modulo $p$ of $P$ lands in $D$.  The fibre $X_P$ is smooth and, generically, its reduction modulo $p$ looks geometrically like the fibre $X_d$.  This doesn't mean that it is non-split, though: the two lines are defined over $\F_p(\sqrt{-a_1/a_2})$, so this variety is split if and only if $-a_1/a_2$ is a square in $\F_p$.  This happens for half of all odd primes $p$.
In this case we find that $X(\Q_p)$ is non-empty only half the time, depending on the prime $p$.
\end{example}

The message from these two contrasting examples is that the geometry of the fibre $X_d$ at the generic point of a divisor $D$ in the base, and in particular whether it is split, strongly affects the local solubility of fibres $X_P$ for $P \in \P^n(k)$ that land in $D$ modulo some prime.
The two obvious questions this raises are: how often does a point $P$ land in a given divisor $D$ modulo some prime?  And do we need to look only at divisors in the base, or also subvarieties of higher codimension?

The answer to the second question is that it turns out to be enough to look at the behaviour above divisors in the base.
The answer to the first question is that, very often, a point $P \in \P^n(k)$ lands in a given divisor $D \subset \P^n$ modulo \emph{some} prime.  After all, in our examples, the only way that $a_0$ can avoid being divisible by any odd prime is for $a_0$ to be plus or minus a power of $2$.  (We also need the extra conditions $p^2 \ndivs a_0$ and $p \ndivs a_1a_2a_3$, but it can be shown that these also hold 100\% of the time.)

\subsection{The results of Loughran and Smeets}

Quantitative results using these ideas relating bad reduction to local solubility have been proved by several authors, most recently Loughran and Smeets~\cite{LS:few}, building on the technique of Serre~\cite{Serre:CRAS-1990}.
To distinguish the two cases seen in our examples, they defined an invariant associated to a family of varieties $\pi \colon X \to \P^n$, as follows.

Firstly, let $D \subset \P^n$ be a prime divisor with generic point $d$.  The fibre $X_d$ is a scheme over the field $k(D)$; let $L/k(D)$ be a finite Galois extension over which all the geometrically irreducible components of $X_d$ are defined.
The Galois group $\Gamma_d = \Gal(L/k(D))$ acts on the set of components.
Define
\[
\delta_D(\pi) = \frac{ \#\{ \gamma \in \Gamma_d \mid \text{$\gamma$ fixes at least one multiplicity-$1$ component of $X_d$} \} } {\#\Gamma_d} .
\]
Note that, if $X_d$ is split, then there is a component fixed by \emph{all} $\gamma \in \Gamma_d$, and therefore we have $\delta_D(\pi)=1$.
However, it is possible to have $\delta_D(\pi)=1$ without $X_d$ being split.

Now define 
\[
\Delta(\pi) = \sum_{D} (1 - \delta_D(\pi) ),
\]
the sum being over all prime divisors of $\P^n$.
This sum is finite: since the generic fibre of $X$ is geometrically integral, there is a non-empty open set $U \subset \P^n$ above which the fibres are geometrically integral, and so $\delta_D(\pi)$ can differ from $1$ only for the finitely many divisors $D$ not contained in $U$.

Armed with this definition, they prove the following.

\begin{theorem}[Loughran, Smeets~\cite{LS:few}]\label{thm:ls1}
Let $\pi \colon X \to \P^n$ be a family of varieties over a number field.  Then we have
\[
N_\mathrm{loc}(B) \ll \frac{B^{n+1}}{(\log B)^{\Delta(\pi)}} \text{ as } B \to \infty.
\]
In particular, if $\Delta(\pi)>0$ then $0\%$ of the varieties in the family are everywhere locally soluble.
\end{theorem}

For the family of Example~\ref{eg:conics}, there are three divisors in the base, given by $\{a_i=0\}$ for $i=0,1,2$, above which the fibre is not split.  
In each case there are two components interchanged by the Galois action, giving $\delta_D(\pi)=\frac{1}{2}$.
We therefore have $\Delta(\pi) = \frac{1}{2} + \frac{1}{2} + \frac{1}{2} = \frac{3}{2}$,
and the theorem agrees with that of Serre (Theorem~\ref{thm:serre}).

The proof of Theorem~\ref{thm:ls1} makes use of the large sieve, a powerful tool for bounding the size of a subset of $\Z^n$ given conditions on its reduction modulo each prime.
See~\cite{Serre:LMWT} for a description of the large sieve in this context.

In the case $\Delta=0$, the positive proportion of everywhere locally soluble varieties is indeed a product of local densities:

\begin{theorem}[Loughran, Smeets]
Suppose that $\pi \colon X \to \P^n$ satisfies $\Delta(\pi)=0$ and $X(\ad_k) \neq \emptyset$.  Then the limit
\[
\lim_{B \to \infty} \frac{N_\mathrm{loc}(B)}{N_\mathrm{tot}(B)}
\]
exists, and is equal to the product of the corresponding local factors.
\end{theorem}

This theorem is proved using the sieve of Ekedahl, a sort of infinite version of the Chinese remainder theorem going back to Ekedahl and developed by various authors; see~\cite[Section~3]{BBL:wa} for more information.

\section{The Brauer--Manin obstruction in families}

In the previous section, we used ideas about bad reduction of varieties in families to prove results about local solubility.  In this section, we use similar techniques to prove results about the Brauer--Manin obstruction.
The rough idea of the conclusion is this:  given a family of varieties and a generic element of the Brauer group, then on most fibres it gives an obstruction to weak approximation but no obstruction to the Hasse principle.

\subsection{Bad reduction and the Brauer--Manin obstruction}

The results in this section come from~\cite{Bright:bad}.
We begin with a number field $k$ and a smooth, proper, geometrically irreducible variety $Y$ over $k$.
Fix a place $v$ of $k$, let $\O_v$ be the ring of integers in the completion $k_v$ and let $\F_v$ be the residue field, of characteristic $p$.
Let $Y_v$ be the base change of $Y$ to $k_v$, and choose a model $\Y/\O_v$ of $Y$.
The special fibre $Y_0 = \Y \times_{\O_v} \F_v$ is a variety over $\F_v$, the ``reduction of $Y$ modulo $v$''.
For what follows, we assume that $Y_0$ is \emph{smooth} and geometrically irreducible, but not necessarily proper.  
For example, starting from any model with geometrically integral special fibre, we can form such a model by removing the singular locus of the special fibre.

The special fibre $Y_0$ is a prime divisor on the regular scheme $\Y$, with $Y_v = \Y \setminus Y_0$, so there is a residue map
\[
\partial \colon \Br Y_v \np \to \H^1(Y_0, \Q/\Z) \np,
\]
where $\np$ denotes the prime-to-$p$ part of a torsion Abelian group.
Now let $P \in Y_v(k_v)$ be a point that reduces to a point $P_0 \in Y_0(\F_v)$.  
(Since $\Y$ is not proper over $\O_v$, there may be points of $Y_v(k_v)$ that ``miss'' the special fibre and don't reduce to points of $Y_0(\F_v)$, but we assume that $P$ is not one of those.)
The residue map behaves in a functorial way, and we get a commutative diagram
\begin{center}
\begin{tikzcd}
\Br Y_v \np \arrow[r, "\partial"] \arrow[d, "P"'] & \H^1(Y_0,\Q/\Z) \np \arrow[d, "P_0"'] \\
\Br k_v \arrow[r] & \H^1(\F_v,\Q/\Z) \arrow[r] & \Q/\Z 
\end{tikzcd}
\end{center}
where the maps in the bottom row are isomorphisms and their composition is the invariant map $\inv_v$: see Remark~\ref{inv}.
From this diagram we immediately deduce the following: for $\al \in \Br Y_v$ and $P \in Y_v(k_v)$ reducing to $P_0 \in Y_0(\F_v)$, the evaluation $\al(P)$ depends only on $\partial \al$ and on $P_0$.
The problem of trying to understand the map $Y(k_v) \to \Br k_v$ given by evaluating $\al$ is determined entirely on the special fibre $Y_0$.

\begin{example}\label{ex:tame}
Take $k=\Q$, $v=p$ odd, let $Y$ be a variety over $\Q$ and suppose that the quaternion algebra $\al = (p,f)_{\Q(Y)}$ happens to lie in $\Br Y$.
What is the class $\partial(\al) \in \H^1(Y_0,\Q/\Z)$?  It will have order $2$, so will be represented by some Galois $2$-covering $Z \to Y_0$.

Perhaps unsurprisingly, it turns out that $Z$ is the $2$-covering given by the equation $t^2 = \bar{f}$, where $\bar{f} \in \F_p(Y_0)$ is the reduction modulo $p$ of the function $f$.
(This equation doesn't work very well at points where $\bar{f}$ has a zero or a pole, but the condition that $\al$ lie in $\Br Y$ means that the divisor $(\bar{f})$ is a multiple of $2$, and so we can always multiply $f$ by a square to obtain an equation that works at any given point.)

Take a point $P \in Y_v(k_v)$ reducing to $P_0 \in Y_0(\F_p)$.  Assume $\bar{f}(P_0) \neq 0$.  A standard calculation for the Hilbert symbol (see~\cite[Chapter~3, Theorem~1]{Serre:CA}) gives
\[
\inv_p \al(P) = \inv_p (p, f(P)) = \begin{cases}
0 & \text{if } \bar{f}(P_0) \in (\F_p^\times)^2; \\
\frac{1}{2} & \text{otherwise}.
\end{cases}
\]
So we see that indeed the invariant is determined by $P_0$ and by $\bar{f}$, that is, by $\partial(\al)$.

We can use this explicit description to make a further argument.
The variety $Z \to Y_0$ is geometrically irreducible if and only if $\bar{f}$ does not become a square in $\bar{\F}_p(Y_0)$.
In this case, the Lang--Weil theorem shows that, if $p$ is sufficiently large, then $Z(\F_p)$ is non-empty.  Looking at the equation defining $Z$, this means that there are points $P_0 \in Y(\F_p)$ such that $\bar{f}(P_0)$ is a square in $\F_p$.
On the other hand, let $a \in \F_p^\times$ be a non-square; then the ``twisted'' variety $Z'$ defined locally by $at^2 = \bar{f}$ also has $\F_p$-points, so there are also points $P_0 \in Y(\F_p)$ such that $\bar{f}(P_0)$ is a non-square in $\F_p$.
Putting this together shows that the evaluation map $Y_v(k_v) \to \Q/\Z$ given by $P \mapsto \inv_p \al(P)$ takes both values $0$ and $\frac{1}{2}$.

Finally, we can deduce something about the Brauer--Manin obstruction on $Y$.  Let $(P_v) \in Y(\ad_\Q)$ be an adelic point of $Y$, and consider the sum $\sum_v \inv_v \al(P_v)$.
Each term in this sum is either $0$ or $\frac{1}{2}$, and
we have shown that the term at $v=p$ can take either value.
It follows that $\sum_v \inv_v \al(P_v)$ takes both values $0$ and $\frac{1}{2}$ on $Y(\ad_\Q)$.
We therefore have
\[
\emptyset \neq Y(\ad_\Q)^{\Br} \subsetneq Y(\ad_\Q) 
\]
and so there is a Brauer--Manin obstruction to weak approximation of $Y$, but no Brauer--Manin obstruction to the Hasse principle.
\end{example}

In the above example, a very explicit argument allowed us to show that a certain element
of the Brauer group of a variety obstructs weak approximation but not the Hasse principle.
We would like to generalise this argument.

Firstly, note that constant algebras (elements of $\Br_0 Y$) can intervene to complicate arguments like that above.
In the example, we had an element $\al$ of order $2$ in $\Br Y$, and we showed that $\inv_v \al(P_v)$ took both possible values in $\Q/\Z[2]$ for $P_v \in Y(k_v)$.
If we change $\al$ by adding an element of large order in $\Br_0 Y$, then the order of $Y$ becomes large but the image of the evaluation map does not get any larger: rather, it is translated by a constant.
To compensate for this, we can fix a base point $P_v \in Y(k_v)$ and instead consider the map $Q_v \mapsto \inv_v \al(Q_v) - \inv_v \al(P_v)$, which is insensitive to changing $Y$ by a constant algebra.

With this in mind, we make the following definition.

\begin{definition}
Let $Y$ be a variety over $k$, let $v$ be a place of $k$, and fix a point $P_v \in Y(k_v)$.
A class $\al \in \Br Y / \Br_0 Y$ of order $n$ is \emph{prolific} at $v$ if the map
\begin{gather*}
Y(k_v) \to \frac{1}{n}\Z/\Z \\
Q_v \mapsto \inv_v \al(Q_v) - \inv_v \al(P_v)
\end{gather*}
is surjective.
\end{definition}
Note that the definition is independent of the representative $\al$ chosen, and independent of the point $P_v$.

An easy argument similar to that used in the example shows that, if $\al$ is prolific at any place $v$, then $\al$ gives an obstruction to weak approximation on $Y$ but no obstruction to the Hasse principle.

However, an obstruction to the Hasse principle can be given by more than one algebra, so
we would also like to be able to work with any finite subgroup of $\Br Y / \Br_0 Y$.
The following definition generalises the previous one.

\begin{definition}
Let $Y$ be a variety over $k$, let $v$ be a place of $k$, and fix a point $P_v \in Y(k_v)$.
A finite subgroup $B \subset \Br Y / \Br_0 Y$ is \emph{prolific} at $v$ if the map
\begin{gather*}
Y(k_v) \to \Hom(B,\Q/\Z) \\
Q_v \mapsto (\al \mapsto \inv_v \al(Q_v) - \inv_v \al(P_v))
\end{gather*}
is surjective.
\end{definition}

Another similar argument shows that, if $B$ is prolific, then $B$ gives an obstruction to weak approximation but no obstruction to the Hasse principle.
Sometimes it can be useful to use variations of these definitions working at more than one place of $k$: see~\cite[Section~2]{Bright:hp}.

In Example~\ref{ex:tame} we began with an algebra $(p,f)$ and, under the assumption that $\bar{f}$ was not square over $\bar{\F}_p$, proved that the algebra was prolific at $p$.
A more general argument is as follows.
Going back to the notation of the beginning of this section, let $\bar{\F}_v$ be an algebraic closure of $\F_v$ and let $\bar{Y}_0$ be the base change of $Y_0$ to $\bar{\F}_v$.
Define $\bar{\partial}$ to be the composition
\begin{equation}\label{eq:dbar}
\Br Y_v\np \xrightarrow{\partial} \H^1(Y_0, \Q/\Z) \to \H^1(\bar{Y}_0, \Q/\Z).
\end{equation}
Suppose that $\al \in \Br Y_v$ has order $n$ in $\Br Y_v / \Br_0 Y_v$; then the torsor corresponding to $\partial(\al)$ is geometrically irreducible if and only if $\bar{\partial}$ has order $n$ in $\H^1(\bar{Y}_0, \Q/\Z)$.
If this condition holds, then an argument using Lang--Weil shows that $\al$ is prolific at $v$ assuming that $\F_v$ is large enough.
More generally, if $B \subset \Br Y_v / \Br_0 Y_v$ is such that $\bar{\partial}$ is injective on $B$, and $\F_v$ is large enough, then $B$ is prolific at $v$.
See~\cite[Section~5]{Bright:bad} for precise statements.

\subsection{The Brauer--Manin obstruction in families}

In this section, we combine the results of the previous section with the ideas of Section~\ref{sec:ls} to study how the Brauer--Manin obstruction varies in a family of varieties.
Again, the idea is to relate what happens at singular fibres of the family to the bad reduction of individual smooth varieties in the family.

Let $\pi \colon X \to \P^n$ be a family of varieties over a number field $k$, as in Section~\ref{sec:ls}.
Let $K$ be the function field of $\P^n$ and let $\eta \colon \Spec K \to \P^n$ be the generic point.
We will start with an element of the generic fibre $\Br X_\eta$ and see what happens when it is specialised to the different $\Br X_P$ for $P \in \P^n(k)$.

First let us think about ramification of an element of $\Br X_\eta$.
Prime divisors on $X$ are of two kinds: \emph{horizontal} prime divisors $Z$, satisfying $\pi(Z)=\P^n$, and \emph{vertical} prime divisors $Z$, for which $D=\pi(Z)$ is a prime divisor on $\P^n$.
Restricting to the generic fibre gives a bijection between the horizontal prime divisors on $X$ and the prime divisors on $X_\eta$. 
Applying the purity theorem (Theorem~\ref{thm:purity}) to $X_\eta$ shows the following: an element $\al \in \Br k(X)$ lies in $\Br X_\eta$ if and only if $\partial_Z(\al)$ is zero for all horizontal prime divisors $Z$.
Therefore such $\al \in \Br X_\eta$ are ramified only at vertical prime divisors of $X$.

Suppose that $\al \in \Br X_\eta$ is ramified at vertical divisors $Z_i$ $(i=1,\dotsc,r)$ and set $D_i = \pi(Z_i)$.
For $P \in \P^n(k)$ not lying in any of the $D_i$, we can specialise $\al$ to the fibre $X_P$ and obtain an element of $\Br X_P$.

From now on, suppose that the fibres of $\pi$ above codimension-one points of $\P^n$ are geometrically irreducible; in that case, we have $Z_i = \pi^{-1}(D_i)$ in the above notation.

Let $D$ be a prime divisor on $\P^n$, with generic point $d$.  The fibre $Z_d = X_d$ is the generic fibre of the morphism $Z \to D$, and is a variety over the field $k(d)$.
Fix an algebraic closure $\overline{k(d)}$ of $k(d)$ and let $\bar{d} \colon \Spec \overline{k(d)} \to \P^n$ be the corresponding geometric point.
The base change $Z_{\bar{d}}$ is then a variety over the algebraically closed field $\overline{k(d)} = k(\bar{d})$.
As in the previous section, we define a homomorphism $\bar{\partial}_Z$ to be the composition
\[
\Br X_\eta \xrightarrow{\partial_Z} \H^1(k(Z), \Q/\Z) \to \H^1(k(Z_{\bar{d}}), \Q/\Z).
\]
Now let $P \in \P^n(k)$ be a point, and suppose that $P$ lands in the divisor $D$ modulo some place $v$ of $k$.
That is, if we let $P_0$ and $D_0$ be the reductions of $P$ and $D$, respectively, modulo $v$, then $P_0$ lies in $D_0$.
Let $p$ be the residue characteristic of $v$.
As in the previous section, if we assume that everything is sufficiently generic, then the fibre $Y=X_P$ will be a smooth variety over $k$ whose reduction $Y_0$ modulo $v$ looks geometrically like $Z_{\bar{d}}$.
After maybe removing some singular locus of $Y_0$, it can be proved that there is a commutative diagram
\begin{center}
\begin{tikzcd}
\Br X_\eta \np \arrow[r, "\bar{\partial}_Z"] \arrow[d] & \H^1(k(Z_{\bar{d}}), \Q/\Z) \np \arrow[d, dashed] \\
\Br Y \np \arrow[r, "\bar{\partial}"] & \H^1(\bar{Y}_0, \Q/\Z) \np
\end{tikzcd}
\end{center}
in which the bottom row is the homomorphism of~\eqref{eq:dbar}.
(The vertical dashed arrow is the specialisation map, not defined everywhere, but OK on the image of $\Br X_\eta$.)
Let $\al |_Y$ denote the restriction of the algebra $\al$ to $Y$.
Suppose that $\bar{\partial}_Z(\al)$ is represented by an irreducible torsor over $Z_{\bar{d}}$; the torsor representing $\bar{\partial}(\al |_Y)$ is a specialisation of that, so will also be irreducible if $P_0$ lies sufficiently generically in $D_0$.
Then the results of the previous section can be brought in to show that $\al |_Y$ is prolific at $v$ and therefore gives an obstruction to weak approximation on $Y$, but no obstruction to the Hasse principle.

The relevant quantitative question is now as follows.  Fix $\al \in \Br X_\eta$, and fix a prime divisor $D \subset \P^n$ such that $\bar{\partial}_Z(\al)$ is non-zero, where $Z$ is the vertical prime divisor on $X$ lying above $D$.
How many points $P \in \P^n(k)$ have the property that $X_P$ is smooth, and there exists a place $v$ of $k$ such that $P$ lands in $D$ modulo $v$ in a sufficiently generic way?
We would also like the residue characteristic not to divide the order of $\al$.
This is the same counting problem as in the previous section, and again application of the large sieve shows that the number of $P$ of height at most $B$ failing to satisfy this condition is $\ll B^{n+1}/ \log B$.

Assuming that such an $\al$ exists, this argument shows that 100\% of smooth varieties in the family have a Brauer--Manin obstruction to weak approximation.
In general there may be non-zero $\al \in \Br X_\eta$ that satisfy $\bar{\partial}_Z(\al)=0$ for all vertical prime divisors $Z$: for example, this is true if $\al$ lies in the image of $\Br K \to \Br X_\eta$; or if $\al$ extends to $\Br X$; or if the family is constant and $\al$ is a non-zero constant element of the Brauer group of the fibre.
However, for many families of interest these last two possibilities do not occur: for example, it often happens that the total space $X$ is $k$-rational, so that $\Br X$ is trivial.
Imposing hypotheses to rule out such uninteresting $\al$, we obtain the following theorem.

\begin{theorem}[Bright, Browning, Loughran~\cite{BBL:wa}]
Let $k$ be a number field and let 
 $\pi \colon X \to \P^n$ be a flat, surjective $k$-morphism of finite type, 
with $X$ smooth, projective and geometrically integral over $k$.
Let $\eta \colon \Spec K \to \P^n$ denote the generic point
and suppose that the generic fibre $X_\eta$ is geometrically connected.
Assume the following hypotheses:
\begin{enumerate}
\item\label{hyp2} $X(\ad_k)\neq \emptyset$;
\item\label{hyp1} the fibre of $\pi$ at each codimension $1$ point of $\P^n$ is geometrically integral;
\item\label{hyp3} the fibre of $\pi$ at each codimension $2$ point of $\P^n$ has a geometrically reduced component;
\item\label{hyp7} $\H^1(k, \Pic\bar{X})=0$;
\item\label{hyp5} $\Br \bar{X}=0$;
\item\label{hyp6}
$\H^2(k,\Pic \P_{\bar k}^n)\to \H^2(k,\Pic \bar X)$ is injective;
\item\label{hyp4} $\Br X_\eta / \Br K \neq 0$.
\end{enumerate}
Then, for 100\% of points $P \in \P^n(k)$ such that the fibre $X_P$ is smooth and everywhere locally soluble,
there is a Brauer--Manin obstruction to weak approximation on $X_P$.
\end{theorem}

In fact, the theorem in~\cite{BBL:wa} is rather more general, for the following reason.
We saw in~\eqref{eq:h1pic} that, for a variety $Y$ over a number field $k$, there is an isomorphism $\Br_1 Y / \Br k \to \H^1(k, \Pic \bar{Y})$.
However, for a variety over an arbitrary field, such as our $X_\eta$, this homomorphism need not be surjective.
When $X$ is the family of diagonal cubic surfaces mentioned in Theorem~\ref{thm:cubics}, Uematsu~\cite{Uematsu:QJM-2014} has shown that $\Br X_\eta / \Br K$ is trivial, whereas $\H^1(K, \Pic X_{\bar{\eta}})$ has order $3$.
In this case we are never going to obtain interesting elements of the Brauer groups of fibres by specialising elements of $\Br X_\eta$, but we can instead specialise elements of $\H^1(K, \Pic X_{\bar{\eta}})$, and in this way we do obtain non-constant elements in the Brauer groups of the fibres.
The theorem of~\cite{BBL:wa} in fact has, in place of condition~(\ref{hyp4}) given above, the condition that either $\Br X_\eta / \Br K$ or $\H^1(K, \Pic X_{\bar{\eta}})$ be non-zero.  The proof is significantly more complicated, since it requires a new definition of residue maps on $\H^1(K, \Pic X_{\bar{\eta}})$.

To finish, we mention a result which extends these methods to study the Brauer--Manin obstruction to the Hasse principle, not just to weak approximation.
Whereas an obstruction to weak approximation on a variety $Y$ can be proved by showing that a single element of $\Br Y$ is prolific at one place, using the same ideas to show the absence of an obstruction to the Hasse principle requires showing that the whole of $\Br Y/\Br k$ is prolific at some place or set of places.

In the context of a family of varieties, we have only been studying elements of the Brauer groups of fibres $X_P$ that arise by specialising elements of $\Br X_\eta$ or $\H^1(K, \Pic X_{\bar{\eta}})$.
To have any hope of proving absence of a Brauer--Manin obstruction on a fibre $X_P$, this specialisation homomorphism needs to be surjective.
One case in which is is true is when the Brauer groups in question are entirely algebraic, that is, $\Br X_\eta = \Br_1 X_\eta$.
In this case, Harari~\cite{Harari:DMJ-1994} has proved that the specialisation map $\H^1(K, \Pic X_{\bar{\eta}}) \to \H^1(k, \Pic \bar{X}_P)$ is an isomorphism for $P$ outside a thin set in $\P^n(k)$, so for 100\% of $P$.  If we assume this hypothesis, the same methods as above allow us to prove the following theorem.

\begin{theorem}[Bright~\cite{Bright:hp}]
Let $k$ be a number field and let $\pi \colon X \to \P^n_k$ be a flat, surjective morphism of finite type, with $X$ smooth, projective and geometrically integral.  Let $\eta \colon \Spec K \to \P^n_k$ denote the generic point and $\bar{\eta} \colon \Spec \bar{K} \to \P^n_k$ a geometric point above $\eta$.  Suppose that the geometric generic fibre $X_{\bar{\eta}}$ is connected and has torsion-free Picard group.  Denote by $\Xbar$ the base change of $X$ to an algebraic closure $\kbar$ of $k$.  Assume the following hypotheses:
\begin{enumerate}
\item\label{locsol} $X(\ad_k) \neq \emptyset$;
\item\label{codim1} the fibre of $\pi$ at each codimension-1 point of $\P^n_k$ is geometrically integral;
\item\label{codim2} the fibre of $\pi$ at each codimension-2 point of $\P^n_k$ has a geometrically reduced component;
\item\label{h1zero} $\H^1(k, \Pic \bar{X})=0$;
\item\label{brzero} $\Br \bar{X}=0$;
\item\label{h2inj} $\H^2(k, \Pic \P^n_{\bar{k}}) \to \H^2(k, \Pic \bar{X})$ is injective.
\item $\Br X_{\bar{\eta}}=0$.
\end{enumerate}
Then, for 100\% of rational points $P \in \P^n(k)$ such that $X_P$ is smooth and everywhere locally soluble, there is no Brauer--Manin obstruction to the Hasse principle on $X_P$.
\end{theorem}

\bibliographystyle{abbrv}
\bibliography{martin}

\end{document}